\documentclass[12pt]{article}
\usepackage{amssymb,amsmath}
\usepackage{cases}
\usepackage{amsfonts}
\usepackage{color,xcolor}
\usepackage[left=2.0cm,right=2.0cm,top=2.0cm,bottom=2.0cm]{geometry}
\usepackage[colorlinks,citecolor=blue,urlcolor=blue]{hyperref}

\newtheorem{theorem}{Theorem}[section]
\newtheorem{corollary}{Corollary}[section]
\newtheorem{lemma}{Lemma}[section]
\newtheorem{proposition}{Proposition}[section]

\newtheorem{definition}{Definition}[section]
\newtheorem{remark}{Remark}[section]

\newcommand{\bal}{\begin{align}}
\newcommand{\bbal}{\begin{align*}}
\newcommand{\beq}{\begin{equation}}
\newcommand{\eeq}{\end{equation}}
\newcommand{\bca}{\begin{cases}}
\newcommand{\eca}{\end{cases}}

\newcommand{\pa}{\partial}
\newcommand{\fr}{\frac}

\newcommand{\De}{\Delta}

\newcommand{\dd}{\mathrm{d}}

\newcommand{\R}{\mathbb{R}}

\newcommand{\les}{\lesssim}

\newcommand{\bi}{\Big}
\newcommand{\g}{\big}
\newcommand{\uu}{\widetilde{u}}
\newcommand{\va}{{\widetilde{\varrho}}}
\begin{document}
\title{Non-uniform dependence for the two-component Camassa-Holm shallow water system}

\author{Jinlu Li$^{1}$, Yanghai Yu$^{2,}$\footnote{E-mail: lijinlu@gnnu.edu.cn; yuyanghai214@sina.com(Corresponding author); mathzwp2010@163.com} and Weipeng Zhu$^{3}$\\
\small $^1$ School of Mathematics and Computer Sciences, Gannan Normal University, Ganzhou 341000, China\\
\small $^2$ School of Mathematics and Statistics, Anhui Normal University, Wuhu 241002, China\\
\small $^3$ School of Mathematics and Information Science, Guangzhou University, Guangzhou 510006, China}

\date{}

\maketitle\noindent{\hrulefill}

{\bf Abstract:}
In this paper, we consider the solution map of the initial value problem to the two-component Camassa-Holm equation on the line. We prove that the solution map of this problem is not uniformly continuous in Sobolev spaces $H^s(\R)\times H^{s-1}(\R)$ for $s>3/2$.


\vskip0mm\noindent{\hrulefill}

\section{Introduction}
In this paper, we are concerned with the Cauchy problem for the following integrable two-component Camassa-Holm (2-CH) shallow water system \cite{Constantin2008,fa,Ivanov,Shabat}
\begin{equation}\label{2-ch}
\begin{cases}
m_t+um_x+2u_xm+\rho\rho_x=0, \; &(t,x)\in \R^+\times\R,\\
m=u-u_{xx},\\
\rho_t+u\rho_x+\rho u_x=0,\\
(u,\rho)(0,x)=(u_0,\rho_0),\; &x\in \R.
\end{cases}
\end{equation}

For $\rho\equiv0$, Eq.\eqref{2-ch} reduces to the famous Camassa-Holm (CH) equation modeling the unidirectional propagation of shallow water waves over a flat bottom, which has attracted much attention as a class of integrable shallow water wave equations in recent twenty years. Concerning the local well-posedness, blow-up criterion, existence of blow-up or global solutions, and non-uniform dependence for the Cauchy problem of (CH) in Sobolev spaces and Besov spaces, we refer to \cite{Bre1,Constantin.Escher2,Constantin.Escher3,Constantin.Escher5,Constantin.Lannes,Constantin.Molinet,Constantin.Strauss,Constantin2,d1,d3,H-K,H-K-M,Guo-Yin,Li-Yin1,li2020,Ro} and the references therein.

When the system is coupled and $\rho$ is a non-constant function, System \eqref{2-ch} is called the two-component Camassa-Holm equations. System \eqref{2-ch} was first proposed by Olver and Rosenau in \cite{Olver}. Constantin--Ivanov \cite{Constantin2008} gave a rigorous justification of the derivation of System \eqref{2-ch}, which is a valid approximation to the governing equations for water waves in the shallow water
regime, where $u(t, x)$ represents the horizontal velocity of the fluid, and $\rho(t, x)$ is related to the
free surface elevation from equilibrium with the following boundary assumptions
\bbal
u\to0\quad\text{and}\quad\rho\to1\quad \text{as}\quad |x|\rightarrow \infty.
\end{align*}
They also investigated conditions for wave-breaking and global small solutions to the system. System \eqref{2-ch} is formally integrable \cite{fa,Shabat} as it can be written as a compatibility condition of two linear systems (Lax pair) with a spectral
parameter (see \cite{Gui2010,Gui2011} and the references therein for more details).

Mathematical properties of the two-component Camassa-Holm have been also studied further in many works. Next, we mainly recall some results which are closely related to our problem. Escher--Lechtenfeld--Yin \cite{Escher2007} investigated local well-posedness for the two-component Camassa-Holm system with initial data $(u_0, \rho_0)\in H^s\times H^{s-1}$ with $s\geq2$ by applying Kato's theory and derived some precise blow-up scenarios for
strong solutions to the system. Guan--Yin \cite{Guan} studied global existence and blow-up phenomena for System \eqref{2-ch} with initial data $(u_0, \rho_0-1)\in H^s\times H^{s-1}$ with $s\geq5/2$. Gui--Liu \cite{Gui2011} established the local well-posedness of System \eqref{2-ch} with initial data $(u_0, \rho_0-1)\in B^s_{p,r}\times B^{s-1}_{p,r}$ with $(p,r)\in[1,\infty]$ and $s>\max\{1+1/p,3/2\}$ (especially in the Sobolev space $H^s\times H^{s-1}$ with $s>3/2$ ) and obtained wave breaking for certain initial profiles. Furthermore, they also proved that the data-to-solution map of this problem is continuous in the weaker space $B^{s'}_{p,r}\times B^{s'-1}_{p,r}$ with $s'<s$.

From the PDE's point of view, it is very interesting that the data-to-solution map whether or not of the Cauchy problem for the 2-CH equation \eqref{2-ch} is uniformly continuous. To the best of our knowledge, the issue has not been investigated up to now. In the present paper, the goal is to show the non-uniform continuity of the solution map of System \eqref{2-ch} in Sobolev spaces $H^s(\R)\times H^{s-1}(\R)$ for $s>3/2$.

Applying the pseudo-differential operator $(1-\pa_x^2)^{-1}$ to the 2-CH equation \eqref{2-ch}, one can rewrite \eqref{2-ch} as a quasi-linear nonlocal evolution system of hyperbolic type
\begin{equation}\label{2-ch1}
\begin{cases}
u_t+uu_x=-\pa_x(1-\pa_x^2)^{-1}\big(u^2+\fr12(u_x)^2+\frac12\rho^2\big), \; &(t,x)\in \R^+\times\R,\\
\rho_t+u\rho_x+\rho u_x=0,\\
(u,\rho)(0,x)=(u_0,\rho_0),\; &x\in \R.
\end{cases}
\end{equation}
Denote $\varrho=\rho-1$, we transform the 2-CH equation \eqref{2-CH} equivalently into the following transport type equation
\begin{equation}\label{2-CH}
\begin{cases}
u_t+\pa_x\Lambda^{-2}\varrho+uu_x=\mathbf{P}(u,\varrho), \; &(t,x)\in \R^+\times\R,\\
\varrho_t+u_x+u\varrho_x+\varrho u_x=0, \\
(u,\varrho)(0,x)=(u_0,\varrho_0),\; &x\in \R,
\end{cases}
\end{equation}
where
\begin{equation}\label{2-CH-0}
\mathbf{P}(u,\varrho)=-\pa_x\Lambda^{-2}\bi(u^2+\fr12(u_x)^2+\frac12\varrho^2\bi)\quad\text{with}\quad\Lambda^{-2}=\g(1-\pa^2_x\g)^{-1}.
\end{equation}
Now let us state our main result of this paper.
\begin{theorem}\label{th-camassa}
Let $s>\fr32$.
The solution map $(u_0,\varrho_0)\mapsto \big(u(t),\varrho(t)\big)=\big(\mathbf{S}_t(u_0),\mathbf{S}_t(\varrho_0)\big)$ of the Cauchy problem \eqref{2-CH}--\eqref{2-CH-0} is not uniformly continuous from any bounded subset in $H^{s}(\R)\times H^{s-1}(\R)$ into $\mathcal{C}([0,T];H^{s}(\R))\times \mathcal{C}([0,T]; H^{s-1}(\R))$. More precisely, there exists two sequences of solutions $(u_{n}^1(t),\varrho_{n}^1(t))$ and $(u_{n}^2(t),\varrho_{n}^2(t))$  such that
\bbal
\|u^1_{0,n},u^2_{0,n}\|_{H^{s}}+\|\varrho^1_{0,n},\varrho^2_{0,n}\|_{H^{s-1}}\lesssim 1
\end{align*}
and
\bbal
\lim_{n\rightarrow \infty}\big(\|u_{0,n}^2-u_{0,n}^1\|_{H^{s}}+\|\varrho_{0,n}^2-\varrho_{0,n}^2\|_{H^{s-1}}\big)= 0
\end{align*}
but for any $t\in[0,T_0]$ with small time $T_0$
\bbal
\liminf_{n\rightarrow \infty}\|u_{n}^2(t)-u_{n}^1(t)\|_{H^{s}}\gtrsim t, \quad \liminf_{n\rightarrow \infty}\|\varrho_{n}^2(t)-\varrho_{n}^1(t)\|_{H^{s-1}}\gtrsim t.
\end{align*}
\end{theorem}
\begin{remark}\label{re1}
Theorem \ref{th-camassa} implies that both the data-to-solution maps of \eqref{2-CH}--\eqref{2-CH-0} are not uniformly continuous in $H^{s}(\R)\times H^{s-1}(\R)$, respectively.
\end{remark}
\begin{remark}\label{re2}
The method we used in proving the Theorem \ref{th-camassa} is different from the CH system \cite{li2020} and can be applied equally well to other related coupled system.
\end{remark}

We outline the main ideas in the proof of Theorem \ref{th-camassa}. At first, we shall analysis the linearized system
\begin{equation}\label{lin}
\begin{cases}
u_t+\pa_x\Lambda^{-2}\varrho=0,\\
\varrho_t+u_x=0,\\
(u,\rho)(0,x)=(u_0,\varrho_0).
\end{cases}
\end{equation}
It should be noted that if $\varrho_0=\Lambda u_0$, we can give the explicit solutions of \eqref{lin}, namely,
\bbal
u=e^{-t\pa_x\Lambda^{-1}}u_{0}\quad \text{and}\quad\varrho=e^{-t\pa_x\Lambda^{-1}}\varrho_{0}.
\end{align*}
Also, Eqs.\eqref{lin} share the linearized part of \eqref{2-CH} and can be indeed regarded as good approximations.
We construct two sequences of approximate initial data $(u_{0,n},\Lambda u_{0,n})$ and $(u_{0,n}+f_n,\Lambda(u_{0,n}+f_n))$, where their distance $f_n$ is arbitrarily close. Meanwhile, we obtain two approximate solutions $(u^{\rm{ap}}_n,\varrho^{\rm{ap}}_n)$ and $(\widetilde{u}^{\rm{ap}}_n,\widetilde{\varrho}^{\rm{ap}}_n)$ which solve \eqref{lin} with the above initial data, respectively. Next, we aim to show that one of the mentioned approximate solutions can approximate the solution maps of \eqref{2-CH} but the other not. For more details, see Propositions \ref{pro1} and \ref{pro2}.
Combining the precious steps, we can conclude that the distance of two solution maps at the initial time is converging to zero, while at any later time it is bounded below by a positive constant, which means the solution map is not uniformly continuous.

The remaining of this paper is organized as follows. In Section 2, we list some notations and known results which will be used in the sequel. In Section 3, we present the local well-posedness result and establish some technical lemmas. In Section 4, we prove our main theorem.
\section{Preliminaries}
{\bf Notations.}\quad The notation $A\les B$ (resp., $A \gtrsim B$) means that there exists a harmless positive constant $c$ such that $A \leq cB$ (resp., $A \geq cB$). Given a Banach space $X$, we denote its norm by $\|\cdot\|_{X}$. We use the simplified notation $\|\mathbf{f}_1,\cdots,\mathbf{f}_n\|_X=\|\mathbf{f}_1\|_X+\cdots+\|\mathbf{f}_n\|_X$ if without confusion.
For all $f\in \mathcal{S}'$, the Fourier transform $\mathcal{F}f$ (also denoted by $\hat{f}$) is defined by
$$
\mathcal{F}f(\xi)=\hat{f}(\xi)=\int_{\R}e^{-ix\xi}f(x)\dd x \quad\text{for any}\; \xi\in\R.
$$
The inverse Fourier transform allows us to recover $f$ from $\hat{u}$:
$$
f(x)=\mathcal{F}^{-1}\hat{f}(x)=\frac{1}{2\pi}\int_{\R}e^{ix\xi}\hat{f}(\xi)\dd\xi.
$$
For $s\in\R$, the operator $\Lambda^s=(1-\pa_x^2)^{\fr{s}2}$ is defined by
$$\widehat{\Lambda^s f}(\xi)=(1+|\xi|^2)^{\fr{s}2}\widehat{f}(\xi).$$
Then the nonhomogeneous Sobolev space is defined by
$$\|f\|^2_{H^s}=\|(1-\pa_x^2)^{\fr{s}2}f\|^2_{L^2}=\int_{\R}(1+|\xi|^2)^s|\hat{f}(\xi)|^2\dd \xi.$$

Next, we will recall some facts about the Littlewood-Paley decomposition, the nonhomogeneous Besov spaces and their some useful properties (see \cite{B.C.D} for more details).

There exists a couple of smooth functions $(\chi,\varphi)$ valued in $[0,1]$, such that $\chi$ is supported in the ball $\mathcal{B}\triangleq \{\xi\in\mathbb{R}:|\xi|\leq \frac 4 3\}$, and $\varphi$ is supported in the ring $\mathcal{C}\triangleq \{\xi\in\mathbb{R}:\frac 3 4\leq|\xi|\leq \frac 8 3\}$. Moreover,
\begin{eqnarray*}
\chi(\xi)+\sum_{j\geq0}\varphi(2^{-j}\xi)=1 \quad \mbox{ for any } \xi\in \R.
\end{eqnarray*}
It is easy to show that $\varphi\equiv 1$ for $\frac43\leq |\xi|\leq \frac32$.

For every $u\in \mathcal{S'}(\mathbb{R})$, the inhomogeneous dyadic blocks ${\Delta}_j$ are defined as follows
\begin{numcases}{\Delta_ju=}
0, & if $j\leq-2$;\nonumber\\
\chi(D)u=\mathcal{F}^{-1}(\chi \mathcal{F}u), & if $j=-1$;\nonumber\\
\varphi(2^{-j}D)u=\mathcal{F}^{-1}\g(\varphi(2^{-j}\cdot)\mathcal{F}u\g), & if $j\geq0$.\nonumber
\end{numcases}
In the inhomogeneous case, the following Littlewood-Paley decomposition makes sense
$$
u=\sum_{j\geq-1}{\Delta}_ju\quad \text{for any}\;u\in \mathcal{S'}(\mathbb{R}).
$$

\begin{definition}[See \cite{B.C.D}]
Let $s\in\mathbb{R}$ and $(p,r)\in[1, \infty]^2$. The nonhomogeneous Besov space $B^{s}_{p,r}(\R)$ is defined by
\begin{align*}
B^{s}_{p,r}(\R):=\Big\{f\in \mathcal{S}'(\R):\;\|f\|_{B^{s}_{p,r}(\mathbb{R})}<\infty\Big\},
\end{align*}
where
\begin{numcases}{\|f\|_{B^{s}_{p,r}(\mathbb{R})}=}
\left(\sum_{j\geq-1}2^{sjr}\|\Delta_jf\|^r_{L^p(\mathbb{R})}\right)^{\fr1r}, &if $1\leq r<\infty$,\nonumber\\
\sup_{j\geq-1}2^{sj}\|\Delta_jf\|_{L^p(\mathbb{R})}, &if $r=\infty$.\nonumber
\end{numcases}
\end{definition}
\begin{remark}\label{re3}
We should be mentioned that
\begin{itemize}
  \item the the Besov space $B_{2,2}^s(\R)$ coincides with the Sobolev space $H^s(\R)$;
  \item the following embedding will be often used implicity:
  $$B^s_{p,q}(\R)\hookrightarrow B^t_{p,r}(\R)\quad\text{for}\;s>t\quad\text{or}\quad s=t,1\leq q\leq r\leq\infty.$$
\end{itemize}
\end{remark}
Finally, we give some important properties which will be also often used throughout the paper.
\begin{lemma}[\cite{B.C.D}]\label{le1}
Let $s>\frac32$. Then we have
\bbal
&\|uv\|_{B^{s-2}_{2,2}(\R)}\leq C\|u\|_{B^{s-2}_{2,2}(\R)}\|v\|_{B^{s-1}_{2,2}(\R)}.
\end{align*}
In particular, when $\mathbf{Q}(u)=\pa_x\Lambda^{-2}u^2$, we have
\bbal
&\|\mathbf{Q}(u)-\mathbf{Q}(v)\|_{B^{s-1}_{2,2}(\R)}\leq C\|u-v\|_{B^{s-2}_{2,2}(\R)}\|u+v\|_{B^{s-1}_{2,2}(\R)}.
\end{align*}
\end{lemma}

\begin{lemma}[\cite{B.C.D}]\label{le2}
For $s>0$, $B^s_{2,2}(\R)\cap L^\infty(\R)$ is an algebra. Moreover, $B^{\frac{1}{2}}_{2,1}(\R)\hookrightarrow L^\infty(\R)$, and for any $u,v \in B^s_{2,2}(\R)\cap L^\infty(\R)$, we have
\bbal
&\|uv\|_{B^{s}_{2,2}(\R)}\leq C(\|u\|_{B^{s}_{2,2}(\R)}\|v\|_{L^\infty(\R)}+\|v\|_{B^{s}_{2,2}(\R)}\|u\|_{L^\infty(\R)}).
\end{align*}
\end{lemma}

\begin{lemma}[\cite{B.C.D,Li-Yin2}]\label{le-tr}
Let $\sigma>-\frac12$. Assume that $f_0\in B^\sigma_{2,2}(\R)$, $g\in L^1([0,T]; B^\sigma_{2,2}(\R))$ and
\begin{numcases}{\pa_{x}\mathbf{u}\in}
L^1([0,T]; B^{\sigma-1}_{2,2}(\R)), & if $\sigma>\fr32$;\nonumber \\
L^1([0,T]; B^{\sigma}_{2,2}(\R)), & if $\sigma=\fr32$;\nonumber\\
L^1([0,T]; B^{1/2}_{2,\infty}(\R)\cap L^\infty(\R)), & if $\sigma<\fr32$.\nonumber
\end{numcases}
If $f\in L^\infty([0,T]; B^\sigma_{2,2}(\R))\cap \mathcal{C}([0,T]; \mathcal{S}'(\R))$ solves the following linear transport equation:
\begin{equation*}
\quad \partial_t f+\mathbf{u}\pa_xf=g,\quad \; f|_{t=0} =f_0.
\end{equation*}

\begin{enumerate}
\item Then there exists a constant $C=C(\sigma)$ such that the following statement holds
\begin{equation*}
\|f(t)\|_{B^\sigma_{2,2}(\R)}\leq e^{CV(t)} \Big(\|f_0\|_{B^\sigma_{2,2}(\R)}+\int_0^t e^{-CV(\tau)} \|g(\tau)\|_{B^\sigma_{2,2}(\R)}\mathrm{d}\tau\Big),
\end{equation*}
where
\begin{numcases}{V(t)=}
\int_0^t \|\pa_x\mathbf{u}(\tau)\|_{B^{\sigma-1}_{2,2}(\R)}\mathrm{d}\tau,\; &if $\sigma>{3\over 2}$ ;\nonumber\\
\int_0^t \|\pa_x\mathbf{u}(\tau)\|_{B^{\sigma}_{2,2}(\R)}\mathrm{d}\tau, \; &if $\sigma={3\over 2}$;\nonumber\\
\int_0^t \|\pa_x\mathbf{u}(\tau)\|_{B^{1/2}_{2,\infty}(\R)\cap L^\infty(\R)}\mathrm{d}\tau,\; &if $\sigma<\fr32$.\nonumber
\end{numcases}
\item If $\sigma>0$, then there exists a constant $C=C(\sigma)$ such that the following statement holds
\begin{align*}
&\|f(t)\|_{B^\sigma_{2,2}(\R)}\leq \|f_0\|_{B^\sigma_{2,2}(\R)}+\int_0^t\|g(\tau)\|_{B^\sigma_{2,2}(\R)}\mathrm{d}\tau \\& \quad \quad
+C\int^t_0\Big(\|f(\tau)\|_{B^\sigma_{2,2}(\R)}\|\pa_x\mathbf{u}(\tau)\|_{L^\infty(\R)}+\|\pa_x\mathbf{u}(\tau)\|_{B^{\sigma-1}_{2,2}(\R)}\|\pa_x f(\tau)| |_{L^\infty(\R)}\Big)\mathrm{d}\tau.
\end{align*}
\end{enumerate}
\end{lemma}

Let us recall the local well-posedness result for the 2-CH equation in Besov spaces.

\begin{lemma}[\cite{Gui2011}]\label{le3-1}
Assume that $s>\fr32$. For any initial data $(u_0,\varrho_0)$ which belongs to $$B_R=\g\{(\psi,\phi)\in B_{2,2}^s\times B^{s-1}_{2,2}: \|\psi\|_{B^{s}_{2,2}}+\|\phi\|_{B^{s-1}_{2,2}}\leq R\g\}\quad\text{for any}\;R>0.$$ Then there exists some $T=T(R)>0$ such that the 2-CH equation has a unique solution $\big(u(t),\varrho (t)\big)\in \mathcal{C}([0,T];B^s_{2,2})\times \mathcal{C}([0,T];B^{s-1}_{2,2})$. Moreover, we have
\begin{align}\label{s}
\|u(t)\|_{B^{s}_{2,2}}+\|\varrho (t)\|_{B^{s-1}_{2,2}}\leq C\big(\|u_0\|_{B^{s}_{2,2}}+\|\varrho_0\|_{B^{s-1}_{2,2}}\big).
\end{align}
\end{lemma}

By Lemmas \ref{le1}--\ref{le-tr} and Gronwall's inequality, we can easily obtain the following result. Since the procedure is standard, we shall not go into the details.

\begin{corollary}\label{cor1.1}
Let $\gamma\geq s-1$ and $(u(t,x),\varrho(t,x))$ be the solution to \eqref{2-CH}. Under the assumptions of Lemma \ref{le3-1}, then there holds for all $t\in[0,T]$,
\begin{align}\label{z}
\|u(t)\|_{B^\gamma_{2,2}}+\|\varrho(t)\|_{B^{\gamma-1}_{2,2}}\leq C(T,\|u_0\|_{B^s_{2,2}},\|\varrho_0\|_{B^{s-1}_{2,2}})(\|u_0\|_{B^\gamma_{2,2}}+\|\varrho_0\|_{B^{\gamma-1}_{2,2}}).
\end{align}
\end{corollary}

\section{Construction of Approximate Solutions}

To construct the approximate solution sequence, we need to introduce smooth, radial cut-off functions to localize the frequency region.
Let $\hat{\phi}\in \mathcal{C}^\infty_0(\mathbb{R})$ be an even, real-valued and non-negative funtion on $\R$ and satify
\begin{numcases}{\hat{\phi}(x)=}
1, &if $|x|\leq \frac{1}{4}$,\nonumber\\
0, &if $|x|\geq \frac{1}{2}$.\nonumber
\end{numcases}
Next, we establish the following crucial lemma which will be used later on.

\begin{lemma}\label{le3-4} Let $s\in\R$, $\lambda\in\big(\fr43,\fr32\big)$ and $n\gg1$. Define the high frequency function $u_{0,n}$ and the low frequency function $f_n$ by
\bal
&u_{0,n}=2^{-ns}\phi(x)\sin \big(\lambda2^nx\big),\label{u0}\\
&f_n=2^{-n}\phi(x).\label{f}
\end{align}
Then for any $\sigma\in\R$, there exists a positive constant $C,\widetilde{C}_1,\widetilde{C}_2$ such that
\bbal\label{le7}
&\|\widehat{u_{0,n}}\|_{L^1}\leq C2^{-ns},\;\;\|\widehat{\pa_xu_{0,n}}\|_{L^1}+\|\widehat{\Lambda u_{0,n}}\|_{L^1}\leq C2^{-n(s-1)},\\
&\|u_{0,n}\|_{H^\sigma}\leq C2^{n(\sigma-s)},\\
&\|f_n\|_{L^\infty}\leq C2^{-n},\;\;\|f_n\|_{H^\sigma}\leq C2^{-n},\\
&\liminf_{n\rightarrow \infty}\|f_n\pa_xu_{0,n}\|_{H^s}\geq \widetilde{C}_1, \quad \liminf_{n\rightarrow \infty}\|f_n \pa_x\Lambda u_{0,n}\|_{H^{s-1}}\geq \widetilde{C}_2.
\end{align*}
\end{lemma}
{\bf Proof.}\quad Easy computations give that
\bbal
\widehat{u_{0,n}}=2^{-ns-1}i\bi[\hat{\phi}\big(\xi+\lambda2^n\big)-\hat{\phi}\big(\xi-\lambda2^n\big)\bi],
\end{align*}
which implies
\bbal
\mathrm{supp} \ \widehat{u_{0,n}}\subset \mathcal{C}\triangleq \Big\{\xi\in\R: \ \lambda2^n-\fr12\leq |\xi|\leq \lambda2^n+\fr12\Big\},
\end{align*}
then, we deduce
$$\|\widehat{\pa_xu_{0,n}}\|_{L^1}\leq C\int_{\mathcal{C}}|\xi\|\widehat{u_{0,n}}(\xi)|\dd\xi\leq C2^{-n(s-1)}.$$
Notice that
\bbal
\widehat{\Lambda u_{0,n}}&=(1+|\xi|^2)^{\fr{1}2}\widehat{u_{0,n}},
\end{align*}
then we have
\bbal
\|\widehat{\Lambda u_{0,n}}\|_{L^1}\leq C\int_{\mathcal{C}}|\xi\|\widehat{u_{0,n}}(\xi)|\dd\xi\leq C2^{-n(s-1)}.
\end{align*}
By the definitions of $f_n$ and $u_{0,n}$, and using the reverse triangle inequality, we obtain
\bbal
\|f_n\pa_x\Lambda u_{0,n}\|_{B^{s-1}_{2,2}}&=2^{n(s-1)}\|\De_{n}\big(f_n\pa_x\Lambda u_{0,n}\big)\|_{L^2}=2^{n(s-1)}\|f_n\pa_x\Lambda u_{0,n}\|_{L^2}
\\&=2^{-n}\big\|\lambda\phi\Lambda\big(\phi\cos(\lambda2^n\cdot)\big)+2^{-n}\phi\Lambda\big(\phi'\sin(\lambda2^n\cdot)\big)\big\|_{L^2}
\\&\gtrsim 2^{-n}\big\|\phi(x)\Lambda\big(\phi\cos(\lambda2^n\cdot)\big)\big\|_{L^2}-2^{-2n}\|\phi\Lambda\big(\phi'\sin (\lambda2^n\cdot)\big)\big\|_{L^2}\\
&\gtrsim 2^{-n}\big\|\Lambda\big(\phi^2\cos(\lambda2^n\cdot)\big)\Big\|_{L^2}-2^{-n}\big\|\big[\Lambda,\phi\big]\big(\phi\cos(\lambda2^n\cdot)\big)\big\|_{L^2}-2^{-n}\\
&\gtrsim 2^{-n}\big\|\Lambda\big(\phi^2\cos(\lambda2^n\cdot)\big)\big\|_{L^2}-2^{-n},
\end{align*}
where we have denoted the commutator by $\big[\Lambda,\phi\big]\varphi=\Lambda(\phi\varphi)-\phi\Lambda\varphi$ and used the estimate
$$\big\|\big[\Lambda,\phi\big]\big(\phi\cos(\lambda2^n\cdot)\big)\big\|_{L^2}\les \|\phi'\|_{L^\infty}\|\phi\cos(\lambda2^n\cdot)\|_{L^2}+\|\phi\|_{H^1}\|\phi\cos(\lambda2^n\cdot)\|_{L^\infty}\quad\text{see \cite{Kato}}.$$
Notice that
\bbal
\mathrm{supp} \ \widehat{\phi}\subset \Big\{\xi\in\R: \ 0\leq |\xi|\leq \fr12\Big\}\;\Rightarrow\;\mathrm{supp} \ \widehat{\phi^2}\subset \Big\{\xi\in\R: \ 0\leq |\xi|\leq 1\Big\},
\end{align*}
which implies
\bbal
\mathrm{supp}\ \mathcal{F}\Big(\phi^2\cos(\lambda2^n\cdot)\Big)\subset  \Big\{\xi\in\R: \ \lambda2^n-1\leq |\xi|\leq \lambda2^n+1\Big\}.
\end{align*}
Thus, by Bernssstein's inequality, we deduce
\bbal
\bi\|\Lambda\big(\phi^2\cos(\lambda2^n\cdot)\big)\bi\|_{L^2}&\gtrsim2^{n}\|\phi^2\cos(\lambda2^n\cdot)\|_{L^2}.
\end{align*}
Following the same procedure of the Proof of Lemmas 3.2--3.4 in \cite{li2020}, we can prove the remaining inequalities. Here we omit the details.

Before proceeding on, we give two data-to-solution maps for the Cauchy problem \eqref{2-CH}
\bbal
&\big(u_{0,n},\varrho_{0,n}=\Lambda u_{0,n}\big)\mapsto\big(u_{n},\varrho_{n}\big)\triangleq\big(\mathbf{S}_{t}(u_{0,n}),\mathbf{S}_{t}(\varrho_{0,n})\big),\\
&\big(\widetilde{u}_{0,n},\widetilde{\varrho}_{0,n}=\Lambda\widetilde{u}_{0,n}\big)\mapsto\big(\widetilde{u}_{n},\widetilde{\varrho}_{n}\big)
\triangleq\big(\mathbf{S}_{t}(\widetilde{u}_{0,n}),\mathbf{S}_{t}(\widetilde{\varrho}_{0,n})\big),
\end{align*}
where $\widetilde{u}_{0,n}=u_{0,n}+f_n$ and $u_{0,n},f_n$ is defined by \eqref{u0} and \eqref{f}, respectively.

From Corollary \ref{cor1.1} and Lemma \ref{le3-4}, we have for all $\gamma\geq s-1$,
\bbal
&\|u_n\|_{H^\gamma}+\|\varrho_{n}\|_{H^{\gamma-1}}\les\|u_{0,n}\|_{H^\gamma}\les2^{n(\gamma-s)},\\
&\|\widetilde{u}_{n}\|_{H^\gamma}+\|\widetilde{\varrho}_{n}\|_{H^{\gamma-1}}\les\|u_{0,n}+f_n\|_{H^\gamma}\les2^{n(\gamma-s)}.
\end{align*}
Next, we construct two approximate solutions, where one can approximate $\big(u_{n},\varrho_{n}\big)$ and the other one can not approximate $\big(\widetilde{u}_{n},\widetilde{\varrho}_{n}\big)$.
\subsection{First Approximate Solutions}
We construct the first approximate solutions as follows
\bbal
u_{n}^{\rm{ap}}=e^{-t\pa_x\Lambda^{-1}}u_{0,n}\quad \text{and}\quad\varrho_{n}^{\rm{ap}}=e^{-t\pa_x\Lambda^{-1}}\Lambda u_{0,n}.
\end{align*}
By direct calculation, we show that $(u_{n}^{\rm{ap}},\varrho_{n}^{\rm{ap}})$ satisfies the following equation
\begin{equation}\label{app1}
\begin{cases}
\pa_tu_{n}^{\rm{ap}}+\pa_x\Lambda^{-2}\varrho_{n}^{\rm{ap}}=0, \\
\pa_t\varrho_{n}^{\rm{ap}}+\pa_xu_{n}^{\rm{ap}}=0,
\end{cases}
\end{equation}
with the initial data $\big(u_{n}^{\rm{ap}},\varrho_{n}^{\rm{ap}}\big)(t=0,x)=(u_{0,n},\Lambda u_{0,n}).$

By Lemma \ref{le3-4}, we know that the following fact holds for $\gamma\in\R$
\bal\label{im}
\|u_{n}^{\rm{ap}}\|_{H^\gamma}+\|\varrho_{n}^{\rm{ap}}\|_{H^{\gamma-1}}=\|u_{0,n}\|_{H^\gamma}+\|\varrho_{0,n}\|_{H^{\gamma-1}}\leq C2^{n(\gamma-s)}.
\end{align}
Setting the difference between the real and the approximate solutions
$$v=u_n-u_{n}^{\rm{ap}}\quad\text{and}\quad w=\varrho_n-\varrho_{n}^{\rm{ap}},$$
then we deduce from \eqref{2-CH} and \eqref{app1} that
\begin{equation}\label{er1}
\begin{cases}
v_t+\pa_x\Lambda^{-2}w+u_n\pa_xv=-v\pa_xu_{n}^{\rm{ap}}-u_{n}^{\rm{ap}}\pa_xu_{n}^{\rm{ap}}+\mathbf{P}(u_{n}^{\rm{ap}},\varrho_{n}^{\rm{ap}})\\
\quad+\mathbf{P}(u_n,\varrho_n)-\mathbf{P}(u_{n}^{\rm{ap}},\varrho_{n}^{\rm{ap}}), \\
w_t+v_x+u_n\pa_xw=-v\pa_x\varrho_{n}^{\rm{ap}}-w\pa_xu_n-\varrho_{n}^{\rm{ap}}\pa_x v
-\pa_x\big(\varrho_{n}^{\rm{ap}}u_{n}^{\rm{ap}}\big),\\
(v,w)(t=0,x)=(0,0).
\end{cases}
\end{equation}
\subsection{Error Estimates for $(v,w)$}
\begin{proposition}\label{pro1}
Assume that $(u_n,\varrho_n)$ solves System \eqref{2-CH} with initial data $(u_{0,n},\Lambda u_{0,n})$. Under the assumptions of Theorem \ref{th-camassa}, then we have
\bal\label{u-es}
\|u_n-u_{n}^{\rm{ap}}\|_{B^{s}_{2,2}}+\|\varrho_n-\varrho_{n}^{\rm{ap}}\|_{B^{s-1}_{2,2}}\leq C2^{-\frac{n}{2}(s-\fr32)}.
\end{align}
\end{proposition}
{\bf Proof.}\quad
Using Lemma \ref{le-tr} to \eqref{er1} yields
\bal\label{ljl}
\|v\|_{B^{s-1}_{2,2}}+\|w\|_{B^{s-2}_{2,2}}&\lesssim \int^t_0\big(\|v\|_{B^{s-1}_{2,2}}+\|w\|_{B^{s-2}_{2,2}}\big)\big(\|u_{n}\|_{B^{s}_{2,2}}+1\big)\dd \tau\nonumber\\ &\quad+\int^t_0\big(\|\mathbf{P}(u_n,\varrho_n)-\mathbf{P}(u_{n}^{\rm{ap}},\varrho_{n}^{\rm{ap}})\|_{B^{s-1}_{2,2}}
+\|\pa_x\big(\varrho_{n}^{\rm{ap}}u_{n}^{\rm{ap}}\big)\|_{B^{s-2}_{2,2}}\big)\dd \tau\nonumber\\ &\quad
+\int^t_0\|v\pa_xu_{n}^{\rm{ap}},u_{n}^{\rm{ap}}\pa_xu_{n}^{\rm{ap}},\mathbf{P}(u_{n}^{\rm{ap}},\varrho_{n}^{\rm{ap}})\|_{B^{s-1}_{2,2}}\dd \tau\nonumber\\
&\quad +\int^t_0\|v\pa_x\varrho_{n}^{\rm{ap}},w\pa_xu_n,\varrho_{n}^{\rm{ap}}\pa_x v
,\pa_x\big(\varrho_{n}^{\rm{ap}}u_{n}^{\rm{ap}}\big)\|_{B^{s-2}_{2,2}}\dd \tau.
\end{align}
From Lemmas \ref{le1}--\ref{le2}, we estimate the above terms one by one
\bbal
&\|v\pa_xu_{n}^{\rm{ap}}\|_{B^{s-1}_{2,2}}\lesssim \|v\|_{B^{s-1}_{2,2}}\|u_{n}^{\rm{ap}}\|_{B^{s}_{2,2}},\\
&\|u_{n}^{\rm{ap}}\pa_xu_{n}^{\rm{ap}}\|_{B^{s-1}_{2,2}}\lesssim \|u_{n}^{\rm{ap}}\|_{B^{s}_{2,2}}\|u_{n}^{\rm{ap}}\|_{L^\infty}\lesssim 2^{-ns},\\
&\|\mathbf{P}(u_{n}^{\rm{ap}},\varrho_{n}^{\rm{ap}})\|_{B^{s-1}_{2,2}}\lesssim \|u_{n}^{\rm{ap}}\|_{H^{s-\fr12}}\|\pa_xu_{n}^{\rm{ap}}\|_{L^\infty}+\|\varrho_{n}^{\rm{ap}}\|_{H^{s-\fr32}}\|\varrho_{n}^{\rm{ap}}\|_{L^\infty}\lesssim 2^{-(s-\fr12)n},\\
&\|\mathbf{P}(u_n,\varrho_n)-\mathbf{P}(u_{n}^{\rm{ap}},\varrho_{n}^{\rm{ap}})\|_{B^{s-1}_{2,2}}\lesssim \|v\|_{B^{s-1}_{2,2}}\|u_{n}^{\rm{ap}},u_{n}\|_{B^{s}_{2,2}}+\|w\|_{B^{s-2}_{2,2}}\|\varrho_{n}^{\rm{ap}},\varrho_{n}\|_{B^{s-1}_{2,2}}\\
&\qquad\qquad\qquad\quad\quad\quad\quad\quad\quad\quad\lesssim \|v\|_{B^{s-1}_{2,2}}+\|w\|_{B^{s-2}_{2,2}},\\
&\|v\pa_x\varrho_{n}^{\rm{ap}}\|_{B^{s-2}_{2,2}}\lesssim \|v\|_{B^{s-1}_{2,2}}\|\varrho_{n}^{\rm{ap}}\|_{B^{s-1}_{2,2}}\lesssim\|v\|_{B^{s-1}_{2,2}},\\
&\|w\pa_xu_n\|_{B^{s-2}_{2,2}}\lesssim \|w\|_{B^{s-2}_{2,2}}\|u_n\|_{B^{s}_{2,2}}\lesssim\|w\|_{B^{s-2}_{2,2}},\\
&\|\varrho_{n}^{\rm{ap}}\pa_x v\|_{B^{s-2}_{2,2}}\lesssim \|\varrho_{n}^{\rm{ap}}\|_{B^{s-1}_{2,2}}\|v\|_{B^{s-1}_{2,2}}\lesssim\|v\|_{B^{s-1}_{2,2}},\\
&\|\pa_x\big(\varrho_{n}^{\rm{ap}}u_{n}^{\rm{ap}}\big)\|_{B^{s-2}_{2,2}}\lesssim \|u_{n}^{\rm{ap}}\|_{B^{s-1}_{2,2}}\|\varrho_{n}^{\rm{ap}}\|_{L^\infty}+\|\varrho_{n}^{\rm{ap}}\|_{B^{s-1}_{2,2}}\|u_{n}^{\rm{ap}}\|_{L^\infty}\lesssim 2^{-ns}
\end{align*}
Plugging the above estimates into \eqref{ljl}, we obtain
\bbal
\|v\|_{B^{s-1}_{2,2}}+\|w\|_{B^{s-2}_{2,2}}&\lesssim \int^t_0\big(\|v\|_{B^{s-1}_{2,2}}+\|w\|_{B^{s-2}_{2,2}}\big)\dd \tau +2^{-(s-\fr12)n},
\end{align*}
which implies
\bbal
\|v\|_{B^{s-1}_{2,2}}+\|w\|_{B^{s-2}_{2,2}}\lesssim 2^{-(s-\fr12)n}.
\end{align*}
An interpolation argument leads to
\bbal
\|v\|_{B^{s}_{2,2}}+\|w\|_{B^{s-1}_{2,2}}\lesssim
\|v\|^{\fr12}_{B^{s-1}_{2,2}}\|v\|^{\fr12}_{H^{s+1}}+\|w\|^{\fr12}_{B^{s-2}_{2,2}}\|w\|^{\fr12}_{B^{s}_{2,2}}\lesssim 2^{-\frac{n}{2}(s-\fr32)}.
\end{align*}
Thus, we have finished the proof of Proposition \ref{pro1}.

\subsection{Second Approximate Solutions}
We construct the second approximate solutions as follows
\bbal
\uu_{n}^{\rm{ap}}=e^{-t\pa_x\Lambda^{-1}}(u_{0,n}+f_n)\quad \text{and}\quad\va_{n}^{\rm{ap}}=e^{-t\pa_x\Lambda^{-1}}\Lambda(u_{0,n}+f_n).
\end{align*}
Similarly, we know that $(\uu_{n}^{\rm{ap}},\va_{n}^{\rm{ap}})$ satisfies Eq.\eqref{app1} with the initial data $\big(u_{0,n}+f_n,\;\Lambda(u_{0,n}+f_n)\big).$

Introducing the following quantities
 $$\widetilde{v}=\uu_n-\uu_{n}^{\rm{ap}}-t\widetilde{U}_{n}^{\rm{ap}}\quad\text{with}\quad \widetilde{U}_{n}^{\rm{ap}}=-\uu_{n}^{\rm{ap}}\pa_x\uu_{n}^{\rm{ap}}$$
 and
 $$\widetilde{w}=\va_n-\va_{n}^{\rm{ap}}-t\widetilde{V}_{n}^{\rm{ap}}\quad\text{with}\quad \widetilde{V}_{n}^{\rm{ap}}=-\uu_{n}^{\rm{ap}}\pa_x\va_{n}^{\rm{ap}},$$
then we deduce the following error system
\begin{equation}\label{ap2}
\begin{cases}
\widetilde{v}_t+\pa_x\Lambda^{-2}\widetilde{w}+\uu_n\pa_x\widetilde{v}=-\widetilde{v}
\pa_x(\uu_{n}^{\rm{ap}}+t\widetilde{U}_{n}^{\rm{ap}})-t\widetilde{U}_{n}^{\rm{ap}}\pa_x\uu_{n}^{\rm{ap}}-t\uu_{n}^{\rm{ap}}\pa_x\widetilde{U}_{n}^{\rm{ap}}
\\\qquad \qquad\qquad\qquad \qquad  \; -t^2\widetilde{U}_{n}^{\rm{ap}}\pa_x\widetilde{U}_{n}^{\rm{ap}}-t\pa_t\widetilde{U}_{n}^{\rm{ap}}-t\pa_x\Lambda^{-2}\widetilde{V}_{n}^{\rm{ap}}+\mathbf{P}(\uu_n,\va_n), \\
\widetilde{w}_t+\widetilde{v}_x+\uu_n\pa_x\widetilde{w}=-\widetilde{v}
\pa_x(\va_{n}^{\rm{ap}}+t\widetilde{V}_{n}^{\rm{ap}})-t\widetilde{U}_{n}^{\rm{ap}}\pa_x\va_{n}^{\rm{ap}}-t\uu_{n}^{\rm{ap}}\pa_x\widetilde{V}_{n}^{\rm{ap}}
\\\qquad \qquad\qquad\qquad  \;-t^2\widetilde{U}_{n}^{\rm{ap}}\pa_x\widetilde{V}_{n}^{\rm{ap}}-t\pa_t\widetilde{V}_{n}^{\rm{ap}}-t\pa_x\widetilde{U}_{n}^{\rm{ap}}
-\widetilde{\varrho}_n\pa_x\uu_n,\\
(\widetilde{v},\widetilde{w})(t=0,x)=(0,0).
\end{cases}
\end{equation}
\subsection{Error Estimates for $(\tilde{v},\tilde{w})$}
\begin{proposition}\label{pro2}
Assume that $(\uu_n,\va_n)$ solves System \eqref{2-CH} with the initial data $$(\widetilde{u}_{0,n},\va_{0,n})=\big({u}_{0,n}+f_n,\Lambda ({u}_{0,n}+f_n)\big).$$
Under the assumptions of Theorem \ref{th-camassa}, we have
\bal\label{u}
\|\uu_n-\uu_{n}^{\rm{ap}}-t\widetilde{U}_n^{\rm{ap}}\|_{B^s_{2,2}}+\|\va_n-\va_{n}^{\rm{ap}}-t\widetilde{V}_{n}^{\rm{ap}}\|_{B^{s-1}_{2,2}}\les t^2+2^{-n(s-\frac{3}{2})}.
\end{align}
\end{proposition}
{\bf Proof.}\quad Using Lemma \ref{le-tr} to \eqref{ap2}, we obtain
\bal\label{e2}
 \ \|\widetilde{v}\|_{B^{s}_{2,2}}+\|\widetilde{w}\|_{B^{s-1}_{2,2}}&\les \int^t_0\big(\|\widetilde{v}\|_{B^{s}_{2,2}}+\|\widetilde{w}\|_{B^{s-1}_{2,2}}\big)\big(1+\|\widetilde{u}_n\|_{B^{s}_{2,2}}\big)\dd \tau
+\int^t_0\|\widetilde{v}
\pa_x(\uu_{n}^{\rm{ap}}+t\widetilde{U}_{n}^{\rm{ap}})\|_{B^{s}_{2,2}}\dd \tau \nonumber\\
&\quad+\int^t_0\|\widetilde{v}
\pa_x(\va_{n}^{\rm{ap}}+t\widetilde{V}_{n}^{\rm{ap}})\|_{B^{s-1}_{2,2}}\dd \tau+\int^t_0\big(\|\mathbf{P}(\uu_n,\va_n)\|_{B^{s}_{2,2}}+\|\widetilde{\varrho}_n\pa_x\uu_n\|_{B^{s-1}_{2,2}}\big)\dd \tau\nonumber\\
&\quad+\int^t_0\tau\|\widetilde{U}_{n}^{\rm{ap}}\pa_x\uu_{n}^{\rm{ap}},\uu_{n}^{\rm{ap}}\pa_x\widetilde{U}_{n}^{\rm{ap}} ,\widetilde{U}_{n}^{\rm{ap}}\pa_x\widetilde{U}_{n}^{\rm{ap}},\pa_\tau\widetilde{U}_{n}^{\rm{ap}},\pa_x\Lambda^{-2}\widetilde{V}_{n}^{\rm{ap}}\|_{B^{s}_{2,2}}\dd \tau\nonumber\\
&\quad+\int^t_0\tau\|\widetilde{U}_{n}^{\rm{ap}}\pa_x\va_{n}^{\rm{ap}},\uu_{n}^{\rm{ap}}\pa_x\widetilde{V}_{n}^{\rm{ap}}
,\widetilde{U}_{n}^{\rm{ap}}\pa_x\widetilde{V}_{n}^{\rm{ap}},\pa_\tau\widetilde{V}_{n}^{\rm{ap}},\pa_x\widetilde{U}_{n}^{\rm{ap}}\|_{B^{s-1}_{2,2}}\dd \tau.
\end{align}
First, due to Lemmas \ref{le1}-\ref{le2}, it is easy to obtain that
\bbal
&\|\widetilde{U}^{\rm{ap}}_{n}\|_{B^{s-1}_{2,2}}\les\|\uu^{\rm{ap}}_{n}\|_{B^{s-1}_{2,2}}\|\uu^{\rm{ap}}_{n}\|_{B^{s}_{2,2}}\les2^{-n},\\
&\|\widetilde{V}^{\rm{ap}}_{n}\|_{B^{s-2}_{2,2}}\les\|\uu^{\rm{ap}}_{n}\|_{B^{s-1}_{2,2}}\|\va^{\rm{ap}}_{n}\|_{B^{s-1}_{2,2}}\les2^{-n},\\
&\|\widetilde{U}^{\rm{ap}}_{n}\|_{B^{s}_{2,2}}\les\|\uu^{\rm{ap}}_{n}\|_{B^{s-1}_{2,2}}\|\uu^{\rm{ap}}_{n}\|_{B^{s+1}_{2,2}}\les1,\\
&\|\widetilde{V}^{\rm{ap}}_{n}\|_{B^{s-1}_{2,2}}\les\|\uu^{\rm{ap}}_{n}\|_{B^{s-1}_{2,2}}\|\va^{\rm{ap}}_{n}\|_{B^{s}_{2,2}}\les1,\\
&\|\widetilde{U}^{\rm{ap}}_{n}\|_{B^{s+1}_{2,2}}\les\|\uu^{\rm{ap}}_{n}\|_{B^{s-1}_{2,2}}\|\uu^{\rm{ap}}_{n}\|_{B^{s+2}_{2,2}}\les2^n,\\
&\|\widetilde{V}^{\rm{ap}}_{n}\|_{B^{s}_{2,2}}\les\|\uu^{\rm{ap}}_{n}\|_{B^{s}_{2,2}}\|\va^{\rm{ap}}_{n}\|_{B^{s}_{2,2}}+\|\uu^{\rm{ap}}_{n}\|_{B^{s-1}_{2,2}}\|\va^{\rm{ap}}_{n}\|_{B^{s+1}_{2,2}}
\les 2^{n}.
\end{align*}
Furthermore, we deduce
\begin{equation}\label{ljl1}
\begin{cases}
\|\widetilde{v}\pa_x(\uu_{n}^{\rm{ap}}+t\widetilde{U}_{n}^{\rm{ap}})\|_{B^{s}_{2,2}}
\les\|\widetilde{v}\|_{B^{s-1}_{2,2}}\|\uu_{n}^{\rm{ap}},\widetilde{U}_{n}^{\rm{ap}}\|_{B^{s+1}_{2,2}}+\|\widetilde{v}\|_{B^{s}_{2,2}}
\|\uu_{n}^{\rm{ap}},\widetilde{U}_{n}^{\rm{ap}}\|_{B^{s}_{2,2}}\\
~~~~~~~~~~~~~~~~~~~~~~~~~~~\les2^n\|\widetilde{v}\|_{B^{s-1}_{2,2}}+\|\widetilde{v}\|_{B^{s}_{2,2}},\\
\|\widetilde{v}\pa_x(\va_{n}^{\rm{ap}}+t\widetilde{V}_{n}^{\rm{ap}})\|_{B^{s-1}_{2,2}}\les\|\widetilde{v}\|_{B^{s-1}_{2,2}}\|\va_{n}^{\rm{ap}},\widetilde{V}_{n}^{\rm{ap}}\|_{B^{s}_{2,2}}\les2^n\|\widetilde{v}\|_{B^{s-1}_{2,2}},
\end{cases}
\end{equation}
\begin{equation}\label{ljl1-1}
\begin{cases}
\|\widetilde{v}\pa_x(\uu_{n}^{\rm{ap}}+t\widetilde{U}_{n}^{\rm{ap}})\|_{B^{s-1}_{2,2}}
\les\|\widetilde{v}\|_{B^{s-1}_{2,2}}\|\uu_{n}^{\rm{ap}},\widetilde{U}_{n}^{\rm{ap}}\|_{B^{s}_{2,2}}
\les \|\widetilde{v}\|_{B^{s-1}_{2,2}},\\
\|\widetilde{v}\pa_x(\va_{n}^{\rm{ap}}+t\widetilde{V}_{n}^{\rm{ap}})\|_{B^{s-2}_{2,2}}\les\|\widetilde{v}\|_{B^{s-1}_{2,2}}\|\va_{n}^{\rm{ap}},\widetilde{V}_{n}^{\rm{ap}}\|_{B^{s-1}_{2,2}}\les \|\widetilde{v}\|_{B^{s-1}_{2,2}},
\end{cases}
\end{equation}
\begin{equation}\label{yyh1}
\begin{cases}
\|\widetilde{U}_{n}^{\rm{ap}}\pa_x\uu_{n}^{\rm{ap}}\|_{B^{s-1}_{2,2}}\les2^{-n},\;\|\uu_{n}^{\rm{ap}}\pa_x\widetilde{U}_{n}^{\rm{ap}}\|_{B^{s-1}_{2,2}}\les2^{-n},\\
\|\widetilde{U}_{n}^{\rm{ap}}\pa_x\widetilde{U}_{n}^{\rm{ap}}\|_{B^{s-1}_{2,2}}\les2^{-n},\;\|\widetilde{U}_{n}^{\rm{ap}}\pa_x\va_{n}^{\rm{ap}}\|_{B^{s-2}_{2,2}}\les2^{-n},\\
\|\uu_{n}^{\rm{ap}}\pa_x\widetilde{V}_{n}^{\rm{ap}}\|_{B^{s-2}_{2,2}}\les2^{-n},\;\|\widetilde{U}_{n}^{\rm{ap}}\pa_x\widetilde{V}_{n}^{\rm{ap}}\|_{B^{s-2}_{2,2}}\les2^{-n},\\
\end{cases}
\end{equation}
\begin{equation}\label{yyh2}
\begin{cases}
\|\widetilde{U}_{n}^{\rm{ap}}\pa_x\uu_{n}^{\rm{ap}}\|_{B^{s}_{2,2}}\les1,\;\|\uu_{n}^{\rm{ap}}\pa_x\widetilde{U}_{n}^{\rm{ap}}\|_{B^{s}_{2,2}}\les1,\\
\|\widetilde{U}_{n}^{\rm{ap}}\pa_x\widetilde{U}_{n}^{\rm{ap}}\|_{B^{s}_{2,2}}\les1,\;\|\widetilde{U}_{n}^{\rm{ap}}\pa_x\va_{n}^{\rm{ap}}\|_{B^{s-1}_{2,2}}\les1,\\
\|\uu_{n}^{\rm{ap}}\pa_x\widetilde{V}_{n}^{\rm{ap}}\|_{B^{s-1}_{2,2}}\les1,\;
\|\widetilde{U}_{n}^{\rm{ap}}\pa_x\widetilde{V}_{n}^{\rm{ap}}\|_{B^{s-1}_{2,2}}\les1.
\end{cases}
\end{equation}
Direct calculation gives that
\bbal
&\pa_t\widetilde{U}_{n}^{\rm{ap}}=(\pa_x\Lambda^{-1}\widetilde{u}_{n}^{\rm{ap}})\pa_x\widetilde{u}_{n}^{\rm{ap}}
+\widetilde{u}_{n}^{\rm{ap}}(\pa^2_x\Lambda^{-1}\widetilde{u}_{n}^{\rm{ap}}),\\
&\pa_t\widetilde{V}_{n}^{\rm{ap}}=(\pa_x\Lambda^{-1}\widetilde{u}_{n}^{\rm{ap}})\pa_x{\va}_{n}^{\rm{ap}}
+\widetilde{u}_{n}^{\rm{ap}}(\pa^2_x\Lambda^{-1}{\va}_{n}^{\rm{ap}}),
\end{align*}
which along with Lemmas \ref{le1}--\ref{le2}  imply
\begin{equation}\label{yyh3}
\begin{cases}
\|\pa_t\widetilde{U}_{n}^{\rm{ap}}\|_{B^{s}_{2,2}}\les\|\widetilde{u}_{n}^{\rm{ap}}\|_{B^{s-1}_{2,2}}\|\widetilde{u}_{n}^{\rm{ap}}\|_{B^{s+1}_{2,2}}+\|\widetilde{u}_{n}^{\rm{ap}}\|^2_{B^{s}_{2,2}}\les1,\\
\|\pa_t\widetilde{V}_{n}^{\rm{ap}}\|_{B^{s-1}_{2,2}}\les\|\widetilde{u}_{n}^{\rm{ap}}\|_{B^{s-1}_{2,2}}\|\va_{n}^{\rm{ap}}\|_{B^{s}_{2,2}}\les1,\\
\|\pa_x\widetilde{U}_{n}^{\rm{ap}}\|_{B^{s-1}_{2,2}}\les\|\widetilde{U}_{n}^{\rm{ap}}\|_{B^{s}_{2,2}}\les1,
\end{cases}
\end{equation}
\begin{equation}\label{yyh4}
\begin{cases}
\|\pa_t\widetilde{U}_{n}^{\rm{ap}}\|_{B^{s-1}_{2,2}}\les\|\widetilde{u}_{n}^{\rm{ap}}\|_{B^{s-1}_{2,2}}\|\widetilde{u}_{n}^{\rm{ap}}\|_{B^{s}_{2,2}}\les2^{-n},\\
\|\pa_t\widetilde{V}_{n}^{\rm{ap}}\|_{B^{s-2}_{2,2}}\les\|\widetilde{u}_{n}^{\rm{ap}}\|_{B^{s-1}_{2,2}}\|\va_{n}^{\rm{ap}}\|_{B^{s-1}_{2,2}}\les2^{-n},\\
\|\pa_x\widetilde{U}_{n}^{\rm{ap}}\|_{B^{s-2}_{2,2}}\les\|\widetilde{U}_{n}^{\rm{ap}}\|_{B^{s-1}_{2,2}}\les2^{-n}.
\end{cases}
\end{equation}
For the term $\mathbf{P}(\uu_n,\va_n)$, we decompose it as
\bbal
\mathbf{P}(\uu_n,\va_n)=\mathcal{A}_1+\mathcal{A}_2+\mathcal{A}_3+\mathcal{A}_4
\end{align*}
where
\bbal\mathcal{A}_1&=-\partial_x\Lambda^{-2}\bi((\uu_{n}^{\rm{ap}})^2+\fr12(\partial_x\uu_{n}^{\rm{ap}})^2+\fr12(\va_{n}^{\rm{ap}})^2\bi),\\
\mathcal{A}_2&=-\partial_x\Lambda^{-2}\bi(\widetilde{v}\uu_n+\fr12\partial_x\widetilde{v}\pa_x\uu_n+\fr12\widetilde{w}\va_n\bi),\\
\mathcal{A}_3&=-\partial_x\Lambda^{-2}\bi(\widetilde{v}(\uu_{n}^{\rm{ap}}+t\widetilde{U}_{n}^{\rm{ap}})+
\fr12\pa_x\widetilde{v}\partial_x(\uu_{n}^{\rm{ap}}+t\widetilde{U}_{n}^{\rm{ap}})
+\fr12\widetilde{w}(\va_{n}^{\rm{ap}}+t\widetilde{V}_{n}^{\rm{ap}})\bi),\\
\mathcal{A}_4&=-t\partial_x\Lambda^{-2}\bi(\widetilde{U}_{n}^{\rm{ap}}(2\uu_{n}^{\rm{ap}}+t\widetilde{U}_{n}^{\rm{ap}})+
\fr12\pa_x\widetilde{U}_{n}^{\rm{ap}}\partial_x(2\uu_{n}^{\rm{ap}}+t\widetilde{U}_{n}^{\rm{ap}})+\fr12\widetilde{V}_{n}^{\rm{ap}}(2\va_{n}^{\rm{ap}}+t\widetilde{V}_{n}^{\rm{ap}})\bi).
\end{align*}
Similarly, we have
\begin{equation}\label{yyh5}
\begin{cases}
\|\mathcal{A}_1\|_{B^{s}_{2,2}}\les2^{-n(s-1)},\\
\|\mathcal{A}_2,\mathcal{A}_3\|_{B^{s}_{2,2}}\les\|\widetilde{v}\|_{B^{s}_{2,2}}+\|\widetilde{w}\|_{B^{s-1}_{2,2}},\\
\|\mathcal{A}_4\|_{B^{s}_{2,2}}\les t,
\end{cases}
\end{equation}
\begin{equation}\label{yyh6}
\begin{cases}
\|\mathcal{A}_1\|_{B^{s-1}_{2,2}}\les2^{-n(s-\frac12)},\\
\|\mathcal{A}_2,\mathcal{A}_3\|_{B^{s-1}_{2,2}}\les\|\widetilde{v}\|_{B^{s-1}_{2,2}}+\|\widetilde{w}\|_{B^{s-2}_{2,2}}, \\
\|\mathcal{A}_4\|_{B^{s-1}_{2,2}}\les t2^{-n}.
\end{cases}
\end{equation}
For the term $\widetilde{\varrho}_n\pa_x\uu_n$, we decompose it as
\bbal
\widetilde{\varrho}_n\pa_x\uu_n&=\underbrace{\widetilde{w}\pa_x{\uu}_n+({\va}_{n}^{\rm{ap}}+t\widetilde{V}^{\rm{ap}}_{n})\pa_x\widetilde{v}}_{\mathcal{A}_5}
+\underbrace{t\big(\widetilde{V}^{\rm{ap}}_{n}\pa_x(\uu_{n}^{\rm{ap}}+t\widetilde{U}_{n}^{\rm{ap}})+\va_{n}^{\rm{ap}}\pa_x\widetilde{U}_{n}^{\rm{ap}}\big)}_{\mathcal{A}_6}
+\underbrace{\va_{n}^{\rm{ap}}\pa_x\uu_{n}^{\rm{ap}}}_{\mathcal{A}_7}.
\end{align*}
Then, we have
\begin{equation}\label{yyh7}
\begin{cases}
\|\mathcal{A}_5\|_{B^{s-1}_{2,2}}\les\|\widetilde{v}\|_{B^{s}_{2,2}}+\|\widetilde{w}\|_{B^{s-1}_{2,2}},\\
\|\mathcal{A}_6\|_{B^{s-1}_{2,2}}\les t,\\
\|\mathcal{A}_7\|_{B^{s-1}_{2,2}}\les \|\widetilde{u}_{n}^{\rm{ap}}\|_{B^{s}_{2,2}}\|\va_{n}^{\rm{ap}}\|_{L^\infty}+\|\va_{n}^{\rm{ap}}\|_{B^{s-1}_{2,2}}\|\pa_x\widetilde{u}_{n}^{\rm{ap}}\|_{L^\infty}\lesssim 2^{-n(s-1)},\\
\end{cases}
\end{equation}
\begin{equation}\label{yyh8}
\begin{cases}
\|\mathcal{A}_5\|_{B^{s-2}_{2,2}}\les\|\widetilde{v}\|_{B^{s-1}_{2,2}}+\|\widetilde{w}\|_{B^{s-2}_{2,2}}, \\
\|\mathcal{A}_6\|_{B^{s-2}_{2,2}}\les t2^{-n},\\
\|\mathcal{A}_7\|_{B^{s-2}_{2,2}}\les\|\widetilde{u}_{n}^{\rm{ap}}\|_{H^{s-\fr12}}\|\va_{n}^{\rm{ap}}\|_{L^\infty}+\|\va_{n}^{\rm{ap}}\|_{H^{s-\fr32}}\|\pa_x\widetilde{u}_{n}^{\rm{ap}}\|_{L^\infty}\lesssim 2^{-(s-\fr12)n}.
\end{cases}
\end{equation}
Gathering all the above estimates together with \eqref{e2} and using Gronwall's inequality yields
\bal\label{e3}
 \ \|\widetilde{v}\|_{B^{s}_{2,2}}+\|\widetilde{w}\|_{B^{s-1}_{2,2}}
&\les \int^t_02^n(\|\widetilde{v}\|_{B^{s-1}_{2,2}}+\|\widetilde{w}\|_{B^{s-2}_{2,2}})\dd \tau+t^2+2^{-n(s-1)}.
\end{align}
Using Lemma \ref{le-tr} to \eqref{ap2} once again and combining the above estimates, we obtain
\bbal
&\|\widetilde{v}\|_{B^{s-1}_{2,2}}+\|\widetilde{w}\|_{B^{s-2}_{2,2}}\lesssim \int^t_0(\|\widetilde{v}\|_{B^{s-1}_{2,2}}+\|\widetilde{w}\|_{B^{s-2}_{2,2}})\dd \tau+t^2 2^{-n} +2^{-(s-\fr12)n}.
\end{align*}
Using Gronwall's inequality yields
\bal\label{e4}
\|\widetilde{v}\|_{B^{s-1}_{2,2}}+\|\widetilde{w}\|_{B^{s-2}_{2,2}}\leq Ct^22^{-n}+2^{-(s-\fr12)n}.
\end{align}
Putting \eqref{e4} into \eqref{e3} implies the desired \eqref{u-es}. The proof of Proposition \ref{pro2} is completed.
\section{Non-uniform Continuous Dependence}
With Propositions \ref{pro1}--\ref{pro2} in hand, we can prove Theorem \ref{th-camassa}.

{\bf Behavior at time $t=0$.}\quad
Obviously, we have
\bbal
\|\widetilde{u}_{0,n}-u_{0,n}\|_{B^{s}_{2,2}}+\|\widetilde{\varrho}_{0,n}-\varrho_{0,n}\|_{B^{s-1}_{2,2}}\leq C\|f_n\|_{B^{s}_{2,2}}\leq C2^{-n},
\end{align*}
which means that
\bbal
\lim_{n\to\infty}\big(\|\widetilde{u}_{0,n}-u_{0,n}\|_{B^{s}_{2,2}}+\|\widetilde{\varrho}_{0,n}-\varrho_{0,n}\|_{B^{s-1}_{2,2}}\big)=0.
\end{align*}
{\bf Behavior at time $t>0$.}\quad Denote $\mathbb{A}\triangleq-\pa_x\Lambda^{-1}$. Notice that
\bbal
&\widetilde{u}_{n}-u_{n}=t\widetilde{U}_n^{\rm{ap}}+\widetilde{v}-v+e^{-t\pa_x\Lambda^{-1}}f_n,\nonumber\\
&\widetilde{\varrho}_{n}-\varrho_{n}=t\widetilde{V}_n^{\rm{ap}}+\widetilde{w}-w+e^{-t\pa_x\Lambda^{-1}}\Lambda f_n,
\end{align*}
and
\bbal
-\widetilde{U}_{n}^{\rm{ap}}
&=f_n\pa_xu_{0,n}+\underbrace{\big(e^{t\mathbb{A}}f_n-f_n\big)\pa_xe^{t\mathbb{A}}u_{0,n}}_{\mathcal{B}_1}
+\underbrace{f_n\pa_x\big(e^{t\mathbb{A}}u_{0,n}-u_{0,n}\big)}_{\mathcal{B}_2}\\
&\quad+\underbrace{e^{t\mathbb{A}}u_{0,n}\pa_xe^{t\mathbb{A}}\big(u_{0,n}+f_n\big)}_{\mathcal{B}_3}
+\underbrace{e^{t\mathbb{A}}f_n\pa_xe^{t\mathbb{A}}f_n}_{\mathcal{B}_4},\\
-\widetilde{V}_{n}^{\rm{ap}}
&=f_n\pa_x\Lambda u_{0,n}+\underbrace{\big(e^{t\mathbb{A}}f_n-f_n\big)\pa_xe^{t\mathbb{A}}\Lambda u_{0,n}}_{\mathcal{B}_5}
+\underbrace{f_n\pa_x\Lambda \big(e^{t\mathbb{A}}u_{0,n}-u_{0,n}\big)}_{\mathcal{B}_6}\\
&\quad+\underbrace{e^{t\mathbb{A}}u_{0,n}\pa_xe^{t\mathbb{A}}\Lambda \big(u_{0,n}+f_n\big)}_{\mathcal{B}_7}
+\underbrace{e^{t\mathbb{A}}f_n\pa_x\Lambda e^{t\mathbb{A}}f_n}_{\mathcal{B}_8},
\end{align*}
using the triangle inequality and Propositions \ref{pro1}--\ref{pro2}, we deduce that
\bal
\|\widetilde{u}_{n}-u_{n}\|_{B^{s}_{2,2}}
\geq&~t\|\widetilde{U}_n^{\rm{ap}}\|_{B^{s}_{2,2}}-\|\widetilde{v},v\|_{B^{s}_{2,2}}-\|f_n\|_{B^{s}_{2,2}}\nonumber\\
\geq&~ t\|f_n\pa_xu_{0,n}\|_{B^{s}_{2,2}}-t\|\mathcal{B}_1,\mathcal{B}_2,\mathcal{B}_3,\mathcal{B}_4\|_{B^{s}_{2,2}}-C2^{-\frac12(s-\frac32) n}-Ct^{2}\nonumber\\
\geq&~ t\|f_n\pa_xu_{0,n}\|_{B^{s}_{2,2}}-C2^{-\frac12(s-\frac32) n}-Ct^{2},\label{yyh}\\
\|\widetilde{\varrho}_{n}-\varrho_{n}\|_{B^{s-1}_{2,2}}
\geq&~t\|\widetilde{V}_n^{\rm{ap}}\|_{B^{s-1}_{2,2}}-\|\widetilde{w},w\|_{B^{s-1}_{2,2}}-\|\Lambda f_n\|_{B^{s-1}_{2,2}}\nonumber\\
\geq&~ t\|f_n\pa_x\Lambda u_{0,n}\|_{B^{s}_{2,2}}-t\|\mathcal{B}_5,\mathcal{B}_6,\mathcal{B}_7,\mathcal{B}_8\|_{B^{s-1}_{2,2}}-C2^{-\frac12(s-\frac32) n}-Ct^{2}\nonumber\\
\geq&~ t\|f_n\pa_x\Lambda u_{0,n}\|_{B^{s}_{2,2}}-C2^{-\frac12(s-\frac32) n}-Ct^{2},\label{yyh-li}
\end{align}
where we have performed the following computations
\bbal
\|\mathcal{B}_1\|_{B^{s}_{2,2}}
&\les\|e^{t\mathbb{A}}f_n-f_n\|_{B^{s}_{2,2}}\|u_{0,n}\|_{B^{s+1}_{2,2}}\les t\|f_n\|_{B^{s}_{2,2}}\|u_{0,n}\|_{B^{s+1}_{2,2}}\les t,\\
\|\mathcal{B}_2\|_{B^{s}_{2,2}}
&\les\|f_n\|_{B^{s}_{2,2}}\|e^{t\mathbb{A}}u_{0,n}-u_{0,n}\|_{B^{s+1}_{2,2}}\les t\|f_n\|_{B^{s}_{2,2}}\|u_{0,n}\|_{B^{s+1}_{2,2}}\les t,\\
\|\mathcal{B}_3\|_{B^{s}_{2,2}}&\les \|e^{t\mathbb{A}}u_{0,n}\|_{L^\infty}\|e^{t\mathbb{A}}\big(u_{0,n}+f_n\big)\|_{B^{s+1}_{2,2}}+\|e^{t\mathbb{A}}u_{0,n}\|_{B^{s}_{2,2}}\|e^{t\mathbb{A}}\pa_x\big(u_{0,n}+f_n\big)\|_{L^\infty}\les 2^{-n(s-1)}+2^{-n},\\
\|\mathcal{B}_4\|_{B^{s}_{2,2}}&\les\|e^{t\mathbb{A}}f_n\|^2_{B^{s+1}_{2,2}}\les\|f_n\|^2_{B^{s+1}_{2,2}} \les 2^{-2n},\\
\|\mathcal{B}_5\|_{B^{s-1}_{2,2}}
&\les\|e^{t\mathbb{A}}f_n-f_n\|_{B^{s-1}_{2,2}}\|u_{0,n}\|_{B^{s+1}_{2,2}}\les t\|f_n\|_{B^{s-1}_{2,2}}\|u_{0,n}\|_{B^{s+1}_{2,2}}\les t,\\
\|\mathcal{B}_6\|_{B^{s-1}_{2,2}}
&\les\|f_n\|_{B^{s-1}_{2,2}}\|e^{t\mathbb{A}}u_{0,n}-u_{0,n}\|_{B^{s+1}_{2,2}}\les t\|f_n\|_{B^{s-1}_{2,2}}\|u_{0,n}\|_{B^{s+1}_{2,2}}\les t,\\
\|\mathcal{B}_7\|_{B^{s-1}_{2,2}}&\les \|e^{t\mathbb{A}}u_{0,n}\|_{L^\infty}\|e^{t\mathbb{A}}\big(u_{0,n}+f_n\big)\|_{B^{s+1}_{2,2}}+\|e^{t\mathbb{A}}u_{0,n}\|_{B^{s-1}_{2,2}}\|e^{t\mathbb{A}}\pa_x\Lambda\big(u_{0,n}+f_n\big)\|_{L^\infty}\les 2^{-n(s-1)}+2^{-n},\\
\|\mathcal{B}_8\|_{B^{s-1}_{2,2}}&\les\|e^{t\mathbb{A}}f_n\|^2_{B^{s+1}_{2,2}}\les\|f_n\|^2_{B^{s+1}_{2,2}} \les 2^{-2n}.
\end{align*}
Hence, it follows from \eqref{yyh}--\eqref{yyh1} and Lemma \ref{le3-4} that
\bbal
\liminf_{n\rightarrow \infty}\|\widetilde{u}_{n}-u_{n}\|_{B^s_{2,2}}\gtrsim t\quad\text{and}\quad
\liminf_{n\rightarrow \infty}\|\widetilde{\varrho}_{n}-\varrho_{n}\|_{B^{s-1}_{2,2}}\gtrsim t,
\end{align*}
for $t$ small enough. This completes the proof of Theorem \ref{th-camassa}.

\vspace*{1em}
\noindent\textbf{Acknowledgements.}  J. Li is supported by the National Natural Science Foundation of China (Grant No.11801090). Y. Yu is supported by the Natural Science Foundation of Anhui Province (No.1908085QA05). W. Zhu is partially supported by the National Natural Science Foundation of China (Grant No.11901092) and Natural Science Foundation of Guangdong Province (No.2017A030310634).
%\vspace*{1em}

\end{document}